\documentclass[a4paper,10pt]{article}
\usepackage{amsmath,amscd,amsfonts,amssymb,epsf,latexsym}

\makeatletter

\long\def\@makefnt#1{\parindent 1em\noindent
            \hb@xt@1.8em{\hss\@textsuperscript{}}#1}
\long\def\@ftntext#1{\insert\footins{%
    \reset@font\footnotesize
    \interlinepenalty\interfootnotelinepenalty
    \splittopskip\footnotesep
    \splitmaxdepth \dp\strutbox \floatingpenalty \@MM
    \hsize\columnwidth \@parboxrestore
    \color@begingroup
      \@makefnt{%
        \rule\z@\footnotesep\ignorespaces#1\@finalstrut\strutbox}%
    \color@endgroup}}%
\def\subjclass#1{%
  \@ftntext{2000 {\itshape Mathematics Subject Classification.}\enspace #1.}}
\def\keywords#1{%
  \@ftntext{{\itshape Key words and phrases.}\enspace #1.}}
\makeatother

\def\A{{\mathbb A}}
\def\B{{\mathbb B}}

\def\C{{\mathbb C}}
\def\D{{\mathbb D}}

\def\X{{\mathbb X}}
\def\Y{{\mathbb Y}}

\def\AB{ {\mathbb A}\moins {\mathbb B}}

\def\codim{\mathop{\rm codim}}

\def\moins{\raise 1pt\hbox{{$\scriptstyle -$}}}
\def\plus{\raise 1pt\hbox{{$\scriptstyle +$}} }

\def\phi{\varphi}

\newtheorem{theorem}{Theorem}
\newtheorem{proposition}[theorem]{Proposition}
\newtheorem{lemma}[theorem]{Lemma}
\newtheorem{corollary}[theorem]{Corollary}
\newtheorem{remark}[theorem]{Remark}
\newtheorem{definition}[theorem]{Definition}

\newtheorem{convention}[theorem]{Convention}
\newtheorem{example}[theorem]{Example}

\def\proof{\noindent{\bf Proof.\ }}

\def\qed{~\hbox{$\Box$}}

\def\Aut{\mathop{\rm Aut}}
\def\Diff{\mathop{\rm Diff}}
\def\codim{\mathop{\rm codim}}

\def\rank{\mathop{\rm rank}}

\begin{document}

\title{\bf Thom polynomials and Schur functions:\\
towards the singularities $A_i(-)$}

\author{Piotr Pragacz\thanks{Research supported by T\"UB\.ITAK (during 
the stay at METU in Ankara), and by the Humboldt Stiftung (during the stay 
at the MPIM in Bonn).}\\
\small Institute of Mathematics of Polish Academy of Sciences\\
\small \'Sniadeckich 8, 00-956 Warszawa, Poland\\
\small P.Pragacz@impan.gov.pl}

\date{(29.03.2006; revised 06.12.2007)}

\subjclass{05E05, 14N10, 57R45}

\keywords{Thom polynomials, singularities, global singularity
theory, classes of degeneracy loci, Schur functions, resultants}

\maketitle

\begin{abstract}
We develop algebro-combinatorial tools for computing the Thom
polynomials for the Morin singularities $A_i(-)$ ($i\ge 0$).
The main tool is the function $F^{(i)}_r$ defined as a
combination of Schur functions with certain numerical
specializations of Schur polynomials as their coefficients.
We show that the Thom polynomial ${\cal T}^{A_i}$ for the singularity $A_i$
(any $i$) associated with maps $({\bf C}^{\bullet},0)
\to ({\bf C}^{\bullet+k},0)$, with any parameter $k\ge 0$,
under the assumption that $\Sigma^j=\emptyset$ for all $j\ge 2$,
is given by $F^{(i)}_{k+1}$.
Equivalently, this says that ``the $1$-part'' of ${\cal T}^{A_i}$
equals $F^{(i)}_{k+1}$. We investigate
2 examples when ${\cal T}^{A_i}$ apart from its $1$-part consists also
of the $2$-part being a single Schur function with some multiplicity.
Our computations combine the characterization of Thom polynomials
via the ``method of restriction equations'' of Rim\'anyi et al.
with the techniques of Schur functions.
\end{abstract}

\section{Introduction}\label{intro}

The global behavior of singularities is governed by their {\it Thom
polynomials} (cf. \cite{T}, \cite{Kl}, \cite{AVGL}, \cite{Ka},
\cite{Rim2}, \cite{Ka4}). Knowing the Thom polynomial of a singularity $\eta$,
denoted ${\cal T}^{\eta}$, one can compute the cohomology class
represented by the $\eta$-points of a map.

In the present paper, following a series of papers by Rim\'anyi et al.
\cite{RS}, \cite{Rim1}, \cite{Rim2}, \cite{FR}, \cite{BFR},
we study the Thom polynomials for the singularities
$A_i$ associated with maps $({\bf C}^{\bullet},0)
\to ({\bf C}^{\bullet+k},0)$  with parameter $k\ge 0$.

The way of obtaining the thought Thom polynomial is through the solution
of a system of linear equations, which is fine when we want to find
one concrete Thom polynomial, say, for a fixed $k$. However, if we want
to find the Thom polynomials for a series of singularities, associated
with maps $({\bf C}^{\bullet},0) \to ({\bf C}^{\bullet+k},0)$
with $k$ as a parameter, we have to solve {\it simultaneously} a countable
family of systems of linear equations.
We do it here for the restriction
equations for the above mentioned singularities.
Instead of using
{\it Chern monomial expansions} (as the authors of previous papers
constantly did), we use {\it Schur function expansions}.
This puts a more transparent structure on computations of Thom
polynomials (cf. also \cite{FK}, \cite{P23}).

Another feature of using the Schur function expansions for Thom
polynomials is that all the coefficients are {\it nonnegative}.
This has been recently proved by A. Weber and the author in \cite{PW} (see
also \cite{PW2}).

\smallskip

To be more precise, we use here (the specializations of)
{\it supersymmetric} Schur functions also called
``Schur functions in difference of alphabets'' together with their
three basic properties: {\it vanishing}, {\it cancellation} and {\it
factorization}, (cf. \cite{S}, \cite{BR}, \cite{LS}, \cite{P2}, \cite {PT},
\cite{M}, \cite{FP}, and \cite{L}).
These functions contain resultants among themselves. Their
geometric role was illuminated, e.g., in the
study of ${\cal P}$-{\it ideals of singularities $\Sigma^i$}
(cf. \cite[end of Sect.~2 and Theorem 11]{P3}) which is based on 
the enumerative geometry of degeneracy loci of \cite{P}.
In fact, in the present paper (and in \cite{PA3}), we use the point of view
of this last paper to some extent.
We know by the Thom-Damon theorem that ${\cal T}^{A_i}$
is a ${\bf Z}$-linear combination of Schur functions in $TX^*\moins f^*(TY^*)$.
Given a positive integer $h$, we shall say that a ${\bf Z}$-linear
combination
$$
\sum_I \alpha_I S_I
$$
is an \ $h$-{\it combination} \ if for any partition $I$ appearing nontrivially
the following condition $(*)_h$ holds\footnote{We say that one partition
{\it is contained} in another if this holds for their Young diagrams
(cf. \cite {L}).}: $I$ contains the rectangle partition
$$
(k+h,\ldots, k+h)
$$
($h$ times), but it does not contain the larger Young diagram
$$
(k+h+1,\ldots, k+h+1)
$$
($h+1$ times).
For example, a $1$-combination consists of Schur functions containing
a single row $(k+1)$ but not containing $(k+2,k+2)$;
a $2$-combination consists of Schur functions containing
$(k+2,k+2)$ but not containing $(k+3,k+3,k+3)$ etc.
(An $h$-combination, with the argument ``$TX^*\moins f^*(TY^*)$'',
is a typical universal polynomial supported on the ($\bullet - h)$th
degeneracy locus of the derivative morphism of the tangent vector bundles.)
Since the singularity $A_i$ is of Thom-Boardman type $\Sigma^1$,
we have by \cite[Theorem 10]{P2} (based on the structure
of the ${\cal P}$-ideal of the singularity $\Sigma^1$)
that all partitions in the Schur expansion of ${\cal T}^{A_i}$
contain a single row $(k+1)$.
For a fixed $h$, let us consider the sum of all Schur functions appearing
nontrivially in ${\cal T}^{A_i}$
(multiplied by their coefficients) corresponding to partitions satisfying
$(*)_h$. This $h$-combination will be called the $h$-{\it part of}
${\cal T}^{A_i}$. Of course, ${\cal T}^{A_i}$ is a sum of its $h$-parts.

\smallskip

The main body of this paper is devoted to study the $1$-part of the Thom
polynomial for the singularities $A_i$ associated with maps
$({\bf C}^{\bullet},0) \to ({\bf C}^{\bullet+k},0)$  with parameter
$k\ge 0$. We introduce, via its Schur function expansion, the basic
functions $F(\A,-)$ and $F^{(i)}$. Using the properties of
these functions (Proposition \ref{FBr} and Corollary \ref{CF}), we show
(Theorem \ref{TFir}) that it gives the Thom polynomial for $A_i$ when
$\Sigma^j=\emptyset$ for all $j\ge 2$. Equivalently, it says that the
$1$-part of the Thom polynomial for a generic singularity $A_i$
is equal to $F^{(i)}_{k+1}$.
For $k=0$, this polynomial was given in \cite{Po} in the Chern
monomial basis.

With the help of $F^{(1)}$ and $F^{(2)}$,
we reprove the formulas of Thom \cite{T} and Ronga \cite{Ro} for $A_1$,
$A_2$ and for any parameter $k\ge 0$.

We give also computations of two Thom polynomials having apart from
their $1$-parts also the nontrivial $2$-parts (consisting of single
Schur functions with certain multiplicities). We first reprove the result
of Gaffney \cite{G} for $A_4$ and $k=0$. This was also done by Rim\'anyi
\cite{Rim1}; our approach uses Schur functions.
Then we do the computations for $A_3$ and  $k=1$; this, in turn,
can be considered as an introduction to the general case $A_3$ (any $k$)
in \cite{PA3}.

In our calculations, we use extensively the functorial $\lambda$-ring
approach to symmetric functions developed mainly in Lascoux's
book \cite{L}.

Main results of the present paper were announced in \cite{P23}.

Inspired by the present article, \cite{P23}, \cite{P3}, and \cite{PA3},
\"Ozer \"Ozt\"urk \cite{O} computed the Thom polynomials 
for $A_4$ and $k=2,3$.

\section{Recollections on Thom polynomials}

Our main reference for this section is \cite{Rim2}.
We start with recalling what we shall mean by a ``singularity''.
Let $k\ge 0$ be a fixed integer. By a {\it singularity}
we shall mean an equivalence class of stable germs $({\bf C}^{\bullet},0)
\to ({\bf C}^{\bullet+k},0)$, where $\bullet\in {\bf N}$, under the
equivalence
generated by right-left equivalence (i.e. analytic reparametrizations
of the source and target) and suspension.

We recall\footnote{This statement is usually called the Thom-Damon
theorem \cite{T}, \cite{D}.}
that the {\it Thom polynomial} ${\cal T}^{\eta}$ of a singularity
$\eta$ is a polynomial in the formal variables $c_1, c_2, \ldots$ that
after the substitution
\begin{equation}
c_i=c_i(f^*TY-TX)=[c(f^*TY)/c(TX)]_i\,,
\end{equation}
for a general map $f:X \to Y$ between complex analytic manifolds,
evaluates the Poincar\'e dual of $[V^{\eta}(f)]$, where $V^{\eta}(f)$
is the cycle carried by the closure of the set
\begin{equation}
\{x\in X : \hbox{the singularity of} \ f \ \hbox{at} \ x \
\hbox{is} \ \eta \}\,.
\end{equation}
By {\it codimension of a singularity} $\eta$, $\codim(\eta)$,
we shall mean $\codim_X(V^{\eta}(f))$ for such an $f$. The concept of
the polynomial ${\cal T}^{\eta}$ comes from Thom's fundamental paper
\cite{T}.
For a detailed discussion of the {\it existence} of Thom polynomials,
see, e.g., \cite{AVGL}. Thom polynomials associated with group actions
were studied by Kazarian in \cite{Ka}, \cite{Ka2}, \cite{Ka4}.

According to Mather's classification, singularities are in one-to-one
correspondence with finite dimensional ${\bf C}$-algebras.
We shall use the following notation:

\medskip

-- $A_i$ \ (of Thom-Boardman type $\Sigma^{1_i}$) will stand for the stable
germs with local algebra ${\bf C}[[x]]/(x^{i+1})$, $i\ge 0$;

\medskip

-- $I_{2,2}$ \ (of Thom-Boardman type $\Sigma^2$) for stable germs with
local algebra ${\bf C}[[x,y]]/(xy, x^2+y^2)$\,;

\medskip

-- $III_{2,2}$ \ (of Thom-Boardman type $\Sigma^2$) for stable germs
with local algebra ${\bf C}[[x,y]]/(xy, x^2, y^2)$ (here $k\ge 1$).

\medskip

In the present article, the computations of Thom polynomials
shall use the method which stems from a sequence
of papers by Rim\'anyi et al. \cite{RS}, \cite{Rim1}, \cite{Rim2},
\cite {FR}, \cite{BFR}.
We sketch briefly this approach, refering the interested reader for more
details to these papers, the main references being the last three
mentioned items.

Let $k\ge 0$ be a fixed integer, and let $\eta: ({\bf C}^{\bullet},0) \to
({\bf C}^{\bullet+k},0)$ be a stable singularity with a prototype
$\kappa: ({\bf C}^n,0) \to ({\bf C}^{n+k},0)$. The {\it maximal compact
subgroup of the right-left symmetry group}
\begin{equation}
\Aut \kappa = \{(\phi,\psi) \in \Diff({\bf C}^n,0) \times
\Diff({\bf C}^{n+k},0) : \psi \circ \kappa \circ \phi^{-1} = \kappa \}
\end{equation}
of $\kappa$ will be denoted by $G_\eta$.
Even if $\Aut \kappa$ is much too large to be a finite dimensional
Lie group, the concept of its maximal compact subgroup (up to conjugacy)
can be defined in a sensible way (cf. \cite{J} and \cite{W}).
In fact, $G_\eta$ can be chosen so that the images of its projections
to the factors $\Diff({\bf C}^n,0)$ and $\Diff({\bf C}^{n+k},0)$ are
linear. Its representations via the projections on the source ${\bf C}^n$
and the target ${\bf C}^{n+k}$ will be denoted by $\lambda_1(\eta)$ and
$\lambda_2(\eta)$.
The vector bundles associated with the universal principal
$G_\eta$-bundle $EG_\eta \to BG_\eta$ using the representations
$\lambda_1(\eta)$ and $\lambda_2(\eta)$ will be called
$E_{\eta}'$ and $E_{\eta}$. The {\it total Chern
class of the singularity} $\eta$ is defined
in $H^{*}(BG_\eta, {\bf Z})$ by
\begin{equation}
c(\eta):=\frac{c(E_{\eta})}{c(E_{\eta}')}\,.
\end{equation}
The {\it Euler class} of $\eta$ is defined in
$H^{2\codim(\eta)}(BG_\eta, {\bf Z})$ by
\begin{equation}
e(\eta):=e(E_{\eta}')\,.
\end{equation}

Sometimes, it will be convenient not to work with
the whole maximal compact subgroup $G_\eta$ but with its suitable subgroup;
this subgroup should be, however, as ``close'' to $G_\eta$ as possible
(cf. \cite{Rim2}, p. 502). We shall denote this subgroup by the
same symbol $G_\eta$.

In the following theorem, we collect information from \cite{Rim2},
Theorem 2.4 and \cite{FR}, Theorem 3.5, needed for the calculations
in the present paper.

\begin{theorem}\label{TEq} \ Suppose, for a singularity $\eta$, that
the Euler classes of all singularities of smaller codimension than
$\codim(\eta)$, are not zero-divisors \footnote{This is the so-called
``Euler condition'' ({\it loc.cit.}).}. Then we have

\noindent
(i) \ if \ $\xi\ne \eta$ \ and \ $\codim(\xi)\le \codim(\eta)$, then
 \ ${\cal T}^{\eta}(c(\xi))=0$;

\noindent
(ii) \ ${\cal T}^{\eta}(c(\eta))=e(\eta)$.

\noindent
This system of equations (taken for all such $\xi$'s) determines the
Thom polynomial ${\cal T}^{\eta}$ in a unique way.
\footnote{To make it precise, we need one more condition
that the number of singularities(=contact orbits) of smaller codimension
is finite: we may assume that $\eta$ is a {\it simple} singularity type,
i.e., there is no moduli adjacent to $\eta$.}
\end{theorem}

To use this method of determining the Thom polynomials for
singularities, one needs their classification, see, e.g., \cite{dPW}.

To effectively use Theorem \ref{TEq}, we need to study the maximal
compact subgroups of singularities. We recall the following recipe
from \cite{Rim2} pp. 505--507. Let $\eta$ be a singularity whose
prototype is $\kappa: ({\bf C}^n,0)\to ({\bf C}^{n+k},0)$. The germ
$\kappa$ is the miniversal unfolding of another germ $\beta:
({\bf C}^m,0)\to ({\bf C}^{m+k},0)$ with $d\beta=0$. The group $G_\eta$
is a subgroup of the maximal compact subgroup of
the algebraic automorphism group of
the local algebra $Q_\eta$
of $\eta$ times the unitary group $U(k\moins d)$, where $d$
is the difference between the minimal number of relations and the number
of generators of $Q_\eta$.
With $\beta$ well chosen, $G_\eta$ acts as right-left symmetry group
on $\beta$ with representations $\mu_1$ and $\mu_2$. The representations
$\lambda_1$ and $\lambda_2$ are
\begin{equation}
\lambda_1=\mu_1\oplus \mu_V \ \ \hbox{and} \ \
\lambda_2=\mu_2\oplus \mu_V\,,
\end{equation}
where $\mu_V$ is the representation of $G_\eta$ on the unfolding space
$V={\bf C}^{n-m}$ given, for $\alpha \in V$ and $(\phi,\psi)\in G_\eta$,
by
\begin{equation}
(\phi,\psi) \ \alpha = \psi \circ \alpha \circ \phi^{-1}\,.
\end{equation}
For example, for the singularity of type $A_i$: $({\bf C}^{\bullet},0) \to
({\bf C}^{\bullet+k},0)$, we have $G_{A_i}=U(1)\times U(k)$ with
\begin{equation}
\mu_1=\rho_1, \ \ \mu_2=\rho_1^{i+1}\oplus \rho_k, \ \
\mu_V=\oplus_{j=2}^i \ \rho_1^j \oplus \oplus_{j=1}^i (\rho_k \otimes
\rho_1^{-j})\,,
\end{equation}
where $\rho_j$ denotes the standard representation of
the unitary group $U(j)$. Hence, we obtain assertion (i) of the following

\begin{proposition}\label{Pce}
(i) \ Let $\eta=A_i$; for any $k$, writing $x$ and $y_1$,\ldots, $y_k$
for the Chern roots of the universal bundles on $BU(1)$ and $BU(k)$,
\begin{equation}
c(A_i)=\frac{1+(i+1)x}{1+x}\prod_{j=1}^k (1+y_j)\,,
\end{equation}
\begin{equation}\label{eA}
e(A_i)= i! \ x^i \ \prod_{j=1}^k (y_j-ix)\cdots (y_j-2x)(y_j-x)\,.
\end{equation}

\smallskip

\noindent
(ii) \ Let $\eta=I_{2,2}$.
Denote by $H$ the extension of $U(1)\times U(1)$ by ${\bf Z}/2{\bf Z}$ 
(``the group generated by multiplication on the coordinates and their 
exchange''). For $k=0$, we have $G_\eta=H$. Hence, for the purpose of our
computations we can use $G_\eta=U(1)\times U(1)$.
Writing $x_1, x_2$ for the Chern roots of the universal bundles
on two copies of $BU(1)$,
\begin{equation}
c(I_{2,2})=\frac{(1+2x_1)(1+2x_2)}{(1+x_1)(1+x_2)}\,.
\end{equation}

\smallskip

\noindent
(iii) \ Let $\eta=III_{2,2}$; for $k=1$,
$G_\eta=U(2)$, and writing $x_1, x_2$ for the Chern roots of the universal 
bundles on $BU(2)$, we have
\begin{equation}
c(III_{2,2})=\frac{(1\plus 2x_1)(1\plus 2x_2)(1\plus x_1\plus x_2)}
{(1\plus x_1)(1\plus x_2)}\,.
\end{equation}
\end{proposition}
(Assertions (ii) and (iii) are obtained, in a standard way, following
the instructions of \cite{Rim2}, Sect.~4. Assertion (ii) is proved in
\cite[pp.~506--507]{Rim2}, whereas assertion (iii) stems from
\cite[p.~65]{BFR}.)

\section{Recollections on Schur functions}

In this section, we collect needed notions related to symmetric
functions. We adopt a functorial point of view of \cite{L}. 
Namely, given a commutative ring, we treat symmetric
functions as operators acting on the ring. We shall give here only a very
brief summary of the corresponding material from our previous paper \cite{P3}.

\medskip

For $m\in {\bf N}$, by ``an alphabet $\A_m$'' we shall mean an alphabet
$\A=(a_1,\ldots,a_m)$ \ (of cardinality $m$); ditto for
$\B_n=(b_1,\ldots,b_n)$, $\Y_k=(y_1,\ldots,y_k)$, and $\X_2=(x_1,x_2)$.

\begin{definition}\label{cf}
Given two alphabets $\A$, $\B$, the {\it complete functions} $S_i(\AB)$
are defined by the generating series (with $z$ an extra variable):

\begin{equation}
\sum S_i(\AB) z^i =\prod_{b\in \B} (1\moins bz)/\prod_{a\in \A}
(1\moins az)\,.
\end{equation}
\end{definition}

\begin{convention} We shall often identify an alphabet
$\A=\{a_1,\ldots,a_m\}$ with the sum $a_1+\cdots +a_m$ and perform usual
algebraic operations on such elements. For example, $\A b$ will
denote the alphabet $(a_1b,\ldots,a_mb)$.
We will give priority to the algebraic notation over the set-theoretic one.
\end{convention}

\begin{definition}\label{sf}
Given a partition\footnote{We identify partitions with their
Young diagrams, as is customary.} $I=(0\le i_1\le i_2\le \ldots\le i_s)\in
{\bf N}^s$, and alphabets $\A$ and $\B$, the {\it Schur function}
$S_I(\A \moins \B)$ is \begin{equation}\label{schur}
S_I(\A \moins \B):= \Bigl|
     S_{i_p+p-q}(\A \moins \B) \Bigr|_{1\leq p,q \le s}  \ .
\end{equation}
\end{definition}
These functions are often called {\it supersymmetric Schur functions}
or {\it Schur functions in difference of alphabets}. Their properties
were studied, among others, in \cite{BR}, \cite{LS}, \cite{P2},
\cite {PT}, \cite{M}, \cite{FP}, and \cite{L}.
From the last item, we borrow increasing ``French'' partitions
and the determinant of the form (\ref{schur}) evaluating a Schur function.
We shall use the the simplified notation $i_1i_2\cdots i_s$ for
a partition $(i_1,\ldots,i_s)$.

We have the following {\it cancellation property}:
\begin{equation}
S_I((\A + \C) - (\B + \C))=S_I(\A-\B)\,.
\end{equation}

We identify partitions with their Young diagrams, as is customary.

We record the following property ({\it loc.cit.}), justifying
the notational
remark from the end of Section 2; for a partition $I$,
\begin{equation}
S_I(\AB)= (-1)^{|I|}S_J(\B \moins \A)=S_J(\B^* \moins \A^*)\,,
\end{equation}
where $J$ is the conjugate partition of $I$ (i.e. the consecutive
rows of the diagram of $J$ are the transposed columns of the
diagram of $I$), and
$\A^*$ denotes the alphabet $\{-a_1,-a_2,\ldots \}$.

Fix two positive integers $m$ and $n$.
We shall say that a partition $I=(0<i_1\le i_2\le \cdots \le i_s)$
{\it is contained in} the $(m,n)$-hook if either $s\le m$, or $s> m$
and $i_{s-m}\le n$.
Pictorially, this means that the Young diagram of $I$ is contained
in the ``tickened" hook:

\smallskip

$$
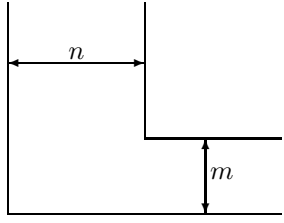
\unitlength=2mm
\begin{picture}(18,14)
\put(0,0){\line(0,1){14}}
\put(0,0){\line(1,0){18}}
\put(9,5){\line(1,0){9}}
\put(9,5){\line(0,1){9}}
\put(4,10){\vector(1,0)5}
\put(4,10){\vector(-1,0)4}
\put(13,2){\vector(0,1)3}
\put(13,2){\vector(0,-1)2}
\put(4.5,10.3){\hbox to0pt{\hss$n$\hss}}
\put(13.3,2.3){\hbox{$m$}}
\end{picture}
$$

\smallskip

We record the following {\it vanishing property}.
Given alphabets $\A$ and $\B$ of cardinalities $m$ and $n$, if
a partition $I$ is not contained in the $(m,n)$-hook,
then ({\it loc.cit.}):
\begin{equation}\label{van}
S_I(\A-\B)=0\,.
\end{equation}

\smallskip

In the present paper, by a {\it symmetric function} we shall mean
a ${\bf Z}$-linear combination of the operators $S_I(-)$.

\medskip
We shall use the following convention from \cite{L1}.

\begin{convention} We may need to specialize a letter to $4$, but this must
not be confused with taking four copies of $1$. To allow one, nevertheless,
specializing a letter to an (integer, or even complex) number $r$ inside
a symmetric function, without introducing intermediate variables,
we write \fbox{$r$} for this specialization. Boxes have to be treated
as single variables. For example, 
$$
S_i(2) = i+1 \ \ \ \hbox{but} \ \ \ S_i(\fbox{2})=2^i\,.
$$
A similar remark applies to ${\bf Z}$-linear combinations of variables.
We have
$$
S_2(\X_2)=x_1^2+ x_1x_2+ x_2^2  \ \ \ \hbox{but} \ \ \
S_2(\fbox{$x_1\plus x_2$})=x_1^2+ 2x_1x_2+ x_2^2\,.
$$
\end{convention}

\begin{definition}
Given two alphabets $\A,\B$, we define their {\it resultant}:
\begin{equation}\label{res}
R(\A,\B):=\prod_{a\in \A,\, b\in \B}(a\moins b)\,.
\end{equation}
\end{definition}
For example, we have the following formal identity:
\begin{equation}\label{fid}
i!(-x)^i\prod_{j=1}^k (ix\moins y_j)\cdots (2x\moins y_j)(x\moins y_j)=
R\bigl(x\plus \fbox{$2x$}\plus \cdots\plus \fbox{$ix$},\Y_k \plus \fbox{$(i\plus 1)x$} \ \bigr)\,.
\end{equation}

We have (cf. \cite{L})
\begin{equation}\label{ER}
R(\A_m,\B_n)= S_{(n^m)}(\AB)=\sum_I S_I(\A) S_{(n^m)/I}(-\B)\,,
\end{equation}
where the sum is over all partitions $I\subset (n^m)$.

When a partition is contained in the $(m,n)$-hook and at the same time it
contains the rectangle $(n^m)$, then we have the following
{\it factorization property} ({\it loc.cit.}):
for partitions $I=(i_1,\ldots,i_m)$ and $J=(j_1,\ldots, j_s)$,
\begin{equation}\label{Fact}
S_{(j_1,\ldots,j_s,i_1+n,\ldots,i_m+n)}(\A_m-\B_n)
=S_I(\A) \ R(\A,\B) \ S_J(-\B)\,.
\end{equation}

\medskip

Rather than the Chern classes
$$
c_i(f^*TY-TX)=[f^*c(TY)/c(TX)]_i \,,
$$
we shall use {\it Segre classes} $S_i$ of the virtual
bundle $TX^*-f^*(TY^*)$,
i.e. complete symmetric functions $S_i(\A -\B)$ for the alphabets
of the {\it Chern roots} $\A, \B$ of $TX^*$ and $TY^*$.

In the present paper, it will be more handy to use, instead of $k$,
a ``shifted'' parameter
\begin{equation}
r:=k+1\,.
\end{equation}
Sometimes, we shall write $\eta(r)$ for the singularity
$\eta: ({\bf C}^{\bullet},0) \to ({\bf C}^{\bullet + r-1},0)$,
and denote the Thom polynomial of $\eta(r)$ by ${\cal T}^{\eta}_r$
-- to emphasize the dependence of both items on $r$.

Note that in our notation, the Thom polynomial for the singularity
$A_1(r)$ for $r\ge 1$, is: ${\cal T}^{A_1}_r=S_r$, instead of $c_{k+1}$
in \cite{Rim2}. In general, a Thom polynomial in terms
of the $c_i$'s (in those papers) will be written here as a linear combination
of Schur functions obtained by changing each $c_i$ to $S_i$ and expanding
in the Schur function basis. Another example is the Thom
polynomial for $A_2(1)$: $c_1^2+c_2$ rewritten in our notation
as ${\cal T}^{A_2}_1=S_{11}+2S_2$.

\medskip

Recall (from the Introduction) that the $h$-{\it part} of ${\cal T}^{A_i}_r$
is the sum of all Schur functions appearing nontrivially
in ${\cal T}^{A_i}_r$ (multiplied by their coefficients) such that the
corresponding partitions satisfy the following condition:
$I$ contains the rectangle partition $\bigl((r\plus h\moins 1)^h\bigr)$,
but it does not contain the larger Young diagram $\bigl((r\plus h)^{h+1}\bigr)$.
The polynomial ${\cal T}^{A_i}_r$ is a sum of its $h$-parts, $h=1,2,\ldots$.

\section{Functions $F(\A, -)$ and $F^{(i)}_r$}

We now pass to the following function $F$ which will give rise to
the $1$-part of ${\cal T}^{A_i}_r$, i.e. to the function $F^{(i)}_r$
that will be studied in this section.
Fix positive integers $m$ and $n$. For an alphabet $\A$ of cardinality $m$,
we define
\begin{equation}
F(\A, -):= \sum_{I} S_I(\A) S_{n-i_m,\ldots,n-i_1,n+|I|}(-)\,,
\end{equation}
where the sum is over partitions $I=(i_1\le i_2 \le \cdots \le i_m \le n)$,
i.e. over $I\subset(n^m)$.

\begin{lemma}\label{LFR} For a variable $x$ and an alphabet $\B$
of cardinality $n$,
\begin{equation}
F(\A,x-\B)= R(x+\A x,\B)\,.
\end{equation}
\end{lemma}
\proof
For a fixed partition $I=(i_1\le i_2 \le \cdots \le i_m\le n)$, it follows
from the factorization property (\ref{Fact}) that
$$
S_{n-i_m,\ldots,n-i_1,n+|I|}(x-\B)=S_{(n^m)/I}(-\B) \ R(x, \B) \ x^{|I|}\,.
$$
Hence, using $S_I(\A x)= S_I(\A) x^{|I|}$, Eq.~(\ref{ER}) and Eq.~(\ref{res}),
we have
$$
\aligned
F(\A,x-\B)=\sum_I S_I(\A) S_{(n^m)/I}(-\B) \ R(x, \B) \ x^{|I|}\\
=\sum_I S_I(\A x) \ S_{(n^m)/I}(-\B) \ R(x, \B)& \\
=R(\A x,\B) \ R(x,\B)
=R(x+\A x, \B)& \,.
\endaligned
$$
The lemma has been proved. 
\qed

\bigskip

The following function $F^{(i)}_r$ will be basic for computing the Thom
polynomials for $A_i$ ($i\ge 1$). We set
\begin{equation}
F^{(i)}_r(-):=\sum_{J} \ S_J(\fbox{$2$}+\fbox{$3$}
+\cdots+\fbox{$i$}) S_{r-j_{i-1},\ldots,r-j_1,r+|J|}(-)\,,
\end{equation}
where the sum is over partitions $J\subset (r^{i-1})$,
and for $i=1$ we understand $F^{(1)}_r(-)=S_r(-)$.

\begin{example} We have
$$
F^{(2)}_r=\sum_{j\le r} S_j(\fbox{$2$}) S_{r-j,r+j}
=\sum_{j\le r} 2^j S_{r-j,r+j}\,;
$$

\smallskip

$$
F^{(3)}_r=\sum_{j_1\le j_2 \le r} S_{j_1,j_2}(\fbox{$2$}+\fbox{$3$})
S_{r-j_2,r-j_1,r+j_1+j_2}\,;
$$
in particular,
$$
F^{(3)}_1=S_{111}+5S_{12}+6S_3
$$
and
$$
F^{(3)}_2=S_{222}+5S_{123}+6S_{114}+19S_{24}+30S_{15}+36S_6;
$$

\smallskip

$$
F^{(4)}_r=\sum_{j_1\le j_2\le j_3 \le r} S_{j_1,j_2,j_3}
(\fbox{$2$}+\fbox{$3$}+\fbox{$4$})
S_{r-j_3,r-j_2,r-j_1,r+j_1+j_2+j_3}\,;
$$
in particular,
$$
F^{(4)}_1=S_{1111}+9S_{112}+26S_{13}+24S_4
$$
and
$$
\aligned
F^{(4)}_2=&S_{2222}+9S_{1223}+26S_{1124}+24S_{1115}+55S_{224}+210S_{125}
+216S_{116}\\
&+391S_{26}+555S_{17}+507S_8\,;
\endaligned
$$

\bigskip

$$
F^{(i)}_1=\sum_{j\le i-1} S_{1^j}(\fbox{$2$}+\fbox{$3$}
+\cdots+\fbox{$i$})S_{1^{i-j-1},j+1}\,.
$$
\end{example}

\medskip

In the following, we shall tacitly assume that $x$, $x_1$, $x_2$,
and $\B_r$ are variables\footnote{Note that these variables will correspond 
in the following to the Chern roots of the {\it cotangent} bundles. On the 
contrary, in Proposition \ref{Pce} the Chern roots of the {\it tangent} bundles
were used. This causes some differences of signs in several formulas.
The same remark applies to our former paper \cite{P3}.} (though many results 
remain valid without this assumption).

The following result gives the key algebraic property of $F^{(i)}_r$.
\begin{proposition}\label{FBr} \ We have
\begin{equation}\label{Br}
F^{(i)}_r(x-\B_r)= R(x+\fbox{$2x$}+\fbox{$3x$}+\cdots
+\fbox{$ix$}\,, \B_r)\,.
\end{equation}
\end{proposition}
\proof
The assertion follows from Lemma \ref{LFR} with $m=i-1$, $n=r$, and
$$
\A=\fbox{$2$}+\fbox{$3$}+\cdots +\fbox{$i$}\,.\qed
$$

\begin{corollary}\label{CF}
Fix an integer $i\ge 1$.

\noindent
(i) For an integer $p\le i$, we have
\begin{equation}\label{eAp}
F^{(i)}_r(x-\B_{r-1}-\fbox{$px$})=0\,.
\end{equation}
(ii) Moreover, we have
\begin{equation}\label{eAn}
F^{(i)}_r(x\moins \B_{r-1}\moins \fbox{$(i\plus 1)x$})=
R(x\plus \fbox{$2x$}\plus \fbox{$3x$}\plus
\cdots\plus\fbox{$ix$}\,, \B_{r-1}\plus \fbox{$(i\plus 1)x$}\, )\,.
\end{equation}
\end{corollary}
\proof Substituting in Eq.~(\ref{Br}):
$$
\B_r=\B_{r-1}+\fbox{$px$}
$$
for $p\le i$, and, respectively,
$$
\B_r=\B_{r-1}+\fbox{$(i\plus 1)x$}\,,
$$
we get the assertions.
\qed

\section{Towards Thom polynomials for $A_i(r)$}

In the following theorem, we shall consider maps $f:X\to Y$ with
degeneracies.

\begin{theorem}\label{TFir}
Suppose that $\Sigma^j(f)=\emptyset$ for $j\ge 2$ \footnote{This says that the
kernel of the derivative map $df: TX \to f^*TY$ of $f$ is a line bundle.}.
Then, for any $r\ge 1$, we have
\begin{equation}
{\cal T}^{A_i}_r = F^{(i)}_r\,.
\end{equation}
\end{theorem}
\proof
By the assumption $\Sigma^j(f)=\emptyset$ for $j\ge 2$, the Euler condition
(needed in Theorem \ref{TEq}) is satisfied here for any $i\ge 0$ and $r\ge 1$.
The equations characterizing ${\cal T}^{A_i}_r$ in the sense of
Theorem \ref{TEq} are, for $p\le i$,
\begin{equation}
P(x-\B_{r-1}-\fbox{$px$})=0\,,
\end{equation}
and additionally, invoking Eq. (\ref{fid}),
\begin{equation}
P(x\moins \B_{r-1}\moins \fbox{$(i\plus 1)x$})=
R(x\plus \fbox{$2x$}\plus \fbox{$3x$}\plus
\cdots\plus\fbox{$ix$}\,, \B_{r-1}\plus \fbox{$(i\plus 1)x$}\, )\,.
\end{equation}
It follows from Corollary \ref{CF} that $P=F^{(i)}_r$ satisfies
these equations. The theorem has been proved.
\qed

\begin{corollary} For any singularity $A_i(r)$, the first part of its
Thom polynomial is equal to $F^{(i)}_r$.
\end{corollary}

\smallskip

\noindent
In the special case $r=1$, Porteous \cite{Po} gave an expression for
the Thom polynomial from the theorem in terms of the Chern monomial basis
(see also \cite{Ku}).

\medskip

The functions $F^{(1)}_r$, $F^{(2)}_r$ give the Thom polynomials for $A_1$, $A_2$
(any $r$) for a general map $f: X \to Y$.

\begin{theorem}(\cite{T}, \cite{Ro}) \
The polynomials \ $S_r$ \ and \ $\sum_{j\le r} 2^j S_{r-j,r+j}$
are Thom polynomials for the singularities $A_1(r)$ and $A_2(r)$.
\end{theorem}
Proof. Since only $A_0$ has smaller codimension than $A_1$,
and only $A_0$, $A_1$ are of smaller codimension
than $A_2$, the Euler conditions hold, and the equations from
Theorem \ref{TEq} characterizing these Thom polynomials are:
\begin{equation}
P(-\B_{r-1})=0, \ \ P(x-\B_{r-1}-\fbox{$2x$})=
R(x, \B_{r-1}\plus \fbox{$2x$})
\end{equation}
for $A_1$, and
\begin{equation}
\aligned
P(-\B_{r-1})=P(x-\B_{r-1}-\fbox{$2x$})=0, \\
P(x\moins \B_{r-1}-\fbox{$3x$})&=
R(x\plus \fbox{$2x$}, \B_{r-1}\plus \fbox{$3x$}\, )
\endaligned
\end{equation}
for $A_2$. Hence the assertion follows from Corollary \ref{CF}.\footnote{Or, 
as the referee points out, it is simpler to say that this follows
from Theorem \ref{TFir} since $\codim {\overline{\Sigma^2}}$ is greater
than $\codim {A_i}$ \ ($i=1,2$).}
\qed

\section{Two examples}

In the present section, we show two (relatively simple) examples of
Schur function expansions of Thom polynomials for $A_i$,
where two $h$-parts appear.
The method used will be applied in \cite{PA3} to more
complicated singularities. Recall that the Thom polynomial ${\cal T}^{A_i}_r$
is a sum of its $h$-parts, the $1$-part being $F^{(i)}_r$.
To get the correct Thom polynomial, one must add to $F^{(i)}_r$
the $h$-parts of ${\cal T}^{A_i}_r$ for $h=2,3,\ldots$.

Let us discuss first $A_4$ for $r=1$ (its codimension is $4$). Then
the singularities $\ne A_4$, whose codimension is $\le \codim(A_4)$
are: $A_0$, $A_1$, $A_2$, $A_3$, $I_{2,2}$. The Thom 
polynomial\footnote{This Thom polynomial was originally computed by Gaffney 
in \cite{G} via the desingularization method. Its alternative derivation 
via solving equations imposed by the above singularities was done by Rim\'anyi 
in \cite{Rim1}). Both authors used Chern monomial expansions.}
is
\begin{equation}\label{A4,1}
{\cal T}^{A_4}_1=S_{1111}+9S_{112}+26S_{13}+24S_4+10S_{22}\,.
\end{equation}
We have
\begin{equation}
F^{(4)}_1=S_{1111}+9S_{112}+26S_{13}+24S_4\,.
\end{equation}
By Corollary \ref{CF}, this function satisfies the following
equations imposed by $A_0$, $A_1$, $A_2$, $A_3$, $A_4$:
\begin{equation}\label{A4}
F^{(4)}_1(0)=F^{(4)}_1(x-\fbox{$2x$})=F^{(4)}_1(x-\fbox{$3x$})=
F^{(4)}_1(x-\fbox{$4x$})=0\,,
\end{equation}

\smallskip

\begin{equation}\label{A4n}
F^{(4)}_1(x-\fbox{$5x$})= R(x+\fbox{$2x$}+\fbox{$3x$}+\fbox{$4x$},
\fbox{$5x$})\,.
\end{equation}
However, $F^{(4)}_1$ does not satisfy the vanishing imposed by $I_{2,2}$.
Namely, we have
\begin{equation}\label{F41}
F^{(4)}_1(\X_2-\fbox{$2x_1$}-\fbox{$2x_2$})
=(-10)x_1x_2(x_1-2x_2)(x_2-2x_1)\,.
\end{equation}\label{F4R}
To see this, invoke Proposition \ref{FBr}:
\begin{equation}
F^{(4)}_1(x-\B_1)=R(x+\fbox{$2x$}+\fbox{$3x$}+\fbox{$4x$},\B_1)\,.
\end{equation}
Substituting to the LHS of Eq.~(\ref{F41}) $x_1=0$, we get by this
proposition
$$
F^{(4)}_1(x_2-\fbox{$2x_2$})=R(x_2+\fbox{$2x_2$}+\fbox{$3x_2$}
+\fbox{$4x_2$},\fbox{$2x_2$})=0\,,
$$
and substituting $x_1=2x_2$,
$$
\aligned
F^{(4)}_1(x_2-\fbox{$2x_1$})=R(x_2+\fbox{$2x_2$}+\fbox{$3x_2$}
+\fbox{$4x_2$},\fbox{$2x_1$})\\
=R(x_2+\fbox{$2x_2$}+\fbox{$3x_2$}+\fbox{$4x_2$},\fbox{$4x_2$})&=0\,.
\endaligned
$$
Therefore 
$$x_1x_2(x_1-2x_2)(x_2-2x_1)
$$
divides this LHS.

To compute the resulting factor we use the specialization $x_1=x_2=1$.
We then have 
$$
x_1x_2(x_1-2x_2)(x_2-2x_1)=1\,,
$$
and $S_{1111}=28$, $S_{112}=-4$, $S_{13}=-1$, $S_4=1$. Hence the factor is
\begin{equation}
F^{(4)}_1=1 \cdot 28 + 9 \cdot 4 + 26 \cdot (-1) + 24 \cdot 1= -10\,,
\end{equation}
and Eq.~(\ref{F41}) is now proved.

On the other hand, the Schur function $S_{22}$ satisfies
Eqs.~(\ref{A4}):
$$
S_{22}(0)=S_{22}(x-\fbox{$2x$})=S_{22}(x-\fbox{$3x$})=S_{22}(x-\fbox{$4x$})
$$
because the partition $22$ is not contained in the $(1,1)$-hook. 
By the same reason, $S_{22}$ satisfies Eq.~(\ref{A4n}) with its RHS 
replaced by zero:
$$
S_{22}(x-\fbox{$5x$})=0\,.
$$
Moreover, we have
\begin{equation}\label{S22}
S_{22}(\X_2-\fbox{$2x_1$}-\fbox{$2x_2$})
=R(\X_2,\fbox{$2x_1$}+\fbox{$2x_2$})
=x_1x_2(x_1-2x_2)(x_2-2x_1)\,.
\end{equation}

Combining Eq.~(\ref{F41}) with Eq.~(\ref{S22}), the desired expression
(\ref{A4,1}) follows.

\bigskip

We now pass to the second example: $A_3$ and $r=2$. The Thom polynomial
in this case was computed originally by Rim\'anyi \cite{Rim2}.
We shall now give its Schur function expansion.
(It is easy to see that the Thom polynomial for $A_3$ and $r=1$ is just equal
to $F^{(3)}_1$.)

Since the singularities $\ne A_3$, whose codimension is $\le \codim(A_3)$
are: $A_0$, $A_1$, $A_2$ and $III_{2,2}$ (cf. \cite{dPW}),
Theorem \ref{TEq} yields the following equations characterizing
${\cal T}^{A_3}_2$, where $b$ is a variable:
\begin{equation}\label{EqA3}
P(-b)=P(x-b-\fbox{$2x$})=P(x-b-\fbox{$3x$})=0\,,
\end{equation}
\begin{equation}\label{A3n}
P(x-b-\fbox{$4x$})=
R(x+\fbox{$2x$}+\fbox{$3x$}, b+\fbox{$4x$}\, )\,
\end{equation}
\begin{equation}
P(\X_2-\D)=0\,.
\end{equation}
Here, 
$$
\D=\fbox{$2x_1$}+\fbox{$2x_2$}+\fbox{$x_1+x_2$}\,.
$$
By Corollary \ref{CF}, the first four equations are satisfied by
the function $F^{(3)}_2$. However $F^{(3)}_2$ does not satisfy the
last vanishing, imposed by $III_{2,2}$.
We shall ``modify'' $F^{(3)}_2$ in order to obtain the Thom polynomial
for $A_3(2)$.

We claim that this Thom polynomial is equal to
\begin{equation}\label{A3(2)}
S_{222}+5S_{123}+6S_{114}+19S_{24}+30S_{15}+36S_6+5S_{33}\,,
\end{equation}
and it differs from its $1$-part $F^{(3)}_2$ by $5S_{33}$ which is
its $2$-part. Indeed, arguing similarly as in the previous example, we have
$$
F^{(3)}_2(\X_2-\D)=(-5)(x_1x_2)^2(x_1-2x_2)(x_2-2x_1)\,.
$$
On the other hand, the Schur function $S_{33}$ satisfies
Eqs.~(\ref{EqA3}):
$$
S_{33}(0)=S_{33}(x-b-\fbox{$2x$})=S_{33}(x-b-\fbox{$3x$})=0
$$
because the partition $33$ is not contained in the $(1,2)$-hook. By the same
reason, $S_{33}$ satisfies Eq.~(\ref{A3n}) with its RHS replaced by zero:
$$
S_{33}(x-b-\fbox{$4x$})=0\,.
$$
Moreover, we have
\begin{equation}
S_{33}(\X_2-\D)
=R(\X_2,\D)
=(x_1x_2)^2(x_1-2x_2)(x_2-2x_1)\,.
\end{equation}

Summing up, we get that the Thom polynomial for $A_3(2)$ has Schur function
expansion (\ref{A3(2)}) indeed.

In \cite{PA3}, we shall give a {\it parametric} Schur function expansion
of the Thom polynomials for the singularities $A_3(r)$ with parameter
$r\ge 1$.

\begin{remark}\rm
Let $\rank({\cal T}^{A_i}_r)$ be the largest $h$ such that there exists
a nontrivial $h$-part in ${\cal T}^{A_i}_r$.
By the results of the present paper, we have

\begin{itemize}

\item $\rank({\cal T}^{A_i}_r)=1$ for $i=1,2$ and any $r$;

\item $\rank({\cal T}^{A_3}_1)=1$, $\rank({\cal T}^{A_3}_2)=2$,
and $\rank({\cal T}^{A_4}_1)=2$.

\end{itemize}

Moreover, we have

\medskip

\begin{itemize}

\item $\rank({\cal T}^{A_3}_r)=2$ for $r\ge 2$ (\cite{P3});

\item $\rank({\cal T}^{A_4}_2)=2$ (\cite{Rim2}, \cite{PW});

\item $\rank({\cal T}^{A_4}_r)=2$ for $r=3,4$ (\cite{O}).

\end{itemize}

Since $\codim(A_i(r))=ir$, for $i\ge 2$ and $r\ge 1$, we clearly have
$$
\rank({\cal T}^{A_i}_r)\le i-1.
$$
This invariant (also for other singularities) will be discussed
in a subsequent paper.
\end{remark}

\smallskip

\noindent
{\bf Acknowledgements} \
Though the author of the present paper is responsible for the exposition
of the details, many computations here were done together with Alain
Lascoux. This help is gratefully acknowledged. The author also thanks
the referee for a careful lecture of the manuscript and pointing out
several corrections. Finally, we note an exceptional relevance of \cite{Jo}
during the work on this paper.

\bigskip

\noindent
{\bf Notes}
1. After the appearance of the first version \cite{P23} of the present paper,
we received a letter from Kazarian \cite{Ka1} informing us that he has
found another formula for ${\cal T}^{A_i}_r$ under the assumptions of
Theorem \ref{TFir}, but modulo a certain ideal (cf. \cite{Ka3}).

\noindent
2. As the referee points out, the Thom polynomials for Morin singularities
have been recently also studied -- using quite different methods --
by Feh\'er and Rim\'anyi in \cite{FR1}, and by B\'erczi and
Szenes in \cite{BSz}.

\end{document}